(10/8/2017, 0h24)

# Inégalité de Turán-Kubilius friable et indépendance asymptotique

Régis de la Bretèche, Youness Lamzouri & Gérald Tenenbaum

**Abstract.** Elaborating on previous works and taking advantage of estimates on the local behaviour of the counting function of friable integers, we determine the optimal range in which the friable Turán-Kubilius constant tends to 1.

**Résumé.** En affinant des résultats antérieurs et en nous appuyant sur des estimations relatives au comportement local de la fonction de comptage des entiers friables, nous déterminons le domaine optimal dans lequel la constante de Turán-Kubilius friable tend vers 1.

## 1. Introduction et énoncé des résultats

Ainsi que mis en évidence par le modèle de Kubilius (cf., par exemple, [7], ch. 3, [8], ch. 12, et [12], ch. III.6.5), les entiers friables, ou sans grand facteur premier, constituent le cadre naturel de l'indépendance asymptotique pour les questions de divisibilité. Un outil essentiel de cette théorie est l'extension friable de l'inégalité de Turán–Kubilius. Abordée dans [1], puis [13], [14], cette étude a été poursuivie et uniformisée dans [3], [11], [9], [4], [5].

Soit $\mathbb{A}$ la classe des fonctions arithmétiques additives complexes. Pour $x \geqslant y \geqslant 2$, nous désignons par $S(x, y)$ l'ensemble des entiers $y$-friables n'excédant pas $x$,[(1)] et par $\alpha = \alpha(x, y)$ le point-selle associé à l'intégrale de Perron pour $\Psi(x, y) := |S(x, y)|$, défini par l'équation transcendante

$$\sum_{p\leqslant y} \frac{\log p}{p^\alpha - 1} = \log x. \tag{1·1}$$

Dans tout ce travail, nous réservons la lettre $p$ pour désigner un nombre premier.

Dans [3], nous avons proposé, pour chaque fonction $f$ de $\mathbb{A}$, un modèle probabiliste $Z_f = Z_{f,x,y}$ de la restriction de $f$ à $S(x, y)$, défini par la formule

$$Z_f := \sum_{p\leqslant y} \xi_p$$

où $\xi_p = \xi_p(f)$ est une variable aléatoire géométrique telle que

$$\mathbb{P}(\xi_p = f(p^\nu)) = g_p(\alpha)/p^{\nu\alpha} \qquad (\nu \in \mathbb{N}), \tag{1·2}$$

les $\xi_p$ étant supposées indépendantes, avec la convention que $f(p^\nu) = 0$ si $p^\nu > x$. Ici et dans la suite, nous utilisons la notation

$$g_m(s) := \prod_{p|m} (1 - 1/p^s) \qquad (m \in \mathbb{N}^*,\ s \in \mathbb{C}).$$

Par convention, si plusieurs valeurs (éventuellement en nombre infini) $f(p^\nu)$ sont égales, la probabilité correspondante figurant au membre de gauche de (1·2) est définie comme la somme des quantités apparaissant au membre de droite.

1. Autrement dit l'ensemble de tous les entiers naturels $n \leqslant x$ dont le plus grand facteur premier est $\leqslant y$.

Comme dans [3], nous définissons la variance semi-empirique $V_f$ de $f$ sur $S(x,y)$ par la formule

$$V_f = V_f(x,y) := \frac{1}{\Psi(x,y)} \sum_{n\in S(x,y)} |f(n) - \mathbb{E}(Z_f)|^2, \tag{1·3}$$

avec donc

$$\mathbb{E}\big(Z_f\big) = \sum_{p^\nu\in S(x,y)} \frac{g_p(\alpha)f(p^\nu)}{p^{\nu\alpha}}\cdot \tag{1·4}$$

Dans son acception la plus forte, l'inégalité de Turán–Kubilius friable consiste en une majoration de $V_f(x,y)$ par un multiple constant de $\mathbb{V}(Z_f)$. Nous avons

$$\mathbb{V}\big(Z_f\big) := \sum_{p^\nu\in S(x,y)} \frac{g_p(\alpha)}{p^{\nu\alpha}}|f(p^\nu)|^2 - \sum_{p\leqslant y}\bigg|\sum_{1\leqslant\nu\leqslant\nu_p} \frac{g_p(\alpha)f(p^\nu)}{p^{\nu\alpha}}\bigg|^2, \tag{1·5}$$

où nous avons posé $\nu_p = \nu_p(x) := \lfloor (\log x)/\log p\rfloor$ $(p \leqslant y)$.

Ainsi qu'il a été observé dans [5], le modèle (1·2) comporte un biais systématique, dû au fait que les probabilités ainsi définies ne sont pas nulles lorsque $p^\nu > x$. Posons

$$w_p = w_p(x,y) := p^{-\alpha\nu_p}, \tag{1·6}$$

de sorte que

$$x^{-\alpha} \leqslant w_p \leqslant \min\{(y/x)^\alpha, x^{-\alpha/2}\}, \qquad (xy)^{-\alpha} \leqslant w_p/p^\alpha \leqslant x^{-\alpha}. \tag{1·7}$$

Nous avons défini dans [5] les variables aléatoires géométriques indépendantes non biaisées $\xi_p^* = \xi_p^*(f)$ vérifiant

$$\mathbb{P}(\xi_p^* = 0) = g_p(\alpha), \qquad \mathbb{P}(\xi_p^* = f(p^\nu)) = \frac{g_p(\alpha)}{(1-w_p)p^{\nu\alpha}} \qquad (1\leqslant\nu\leqslant\nu_p), \tag{1·8}$$

avec la même convention que précédemment pour les cas d'égalité des valeurs prises, et nous posons alors

$$Z_f^* := \sum_{p\leqslant y} \xi_p^*, \tag{1·9}$$

et

$$\begin{aligned} \mathbb{E}(Z_f^*) &= \sum_{p^\nu\in S(x,y)} \frac{g_p(\alpha)f(p^\nu)}{(1-w_p)p^{\nu\alpha}}, \\ \mathbb{V}\big(Z_f^*\big) &:= \sum_{p^\nu\in S(x,y)} \frac{g_p(\alpha)|f(p^\nu)|^2}{(1-w_p)p^{\nu\alpha}} - \sum_{p\leqslant y}\bigg|\sum_{1\leqslant\nu\leqslant\nu_p} \frac{g_p(\alpha)f(p^\nu)}{(1-w_p)p^{\nu\alpha}}\bigg|^2, \end{aligned} \tag{1·10}$$

la variance semi-empirique associée à ce modèle étant définie par

$$V_f^* = V_f^*(x,y) := \frac{1}{\Psi(x,y)} \sum_{n\in S(x,y)} \big|f(n) - \mathbb{E}(Z_f^*)\big|^2. \tag{1·11}$$

Les résultats de [3], [4], et [5] fournissent les majorations universelles

$$(1{\cdot}12)\quad C(x,y):=\sup_{f\in\mathbb{A}} V_f/\mathbb{V}(Z_f)\ll 1,\quad C^*(x,y):=\sup_{f\in\mathbb{A}} V_f^*/\mathbb{V}(Z_f^*)\ll 1\qquad (x\geqslant y\geqslant 2).$$

De plus, il a été établi dans [3] que

$$(1{\cdot}13)\qquad C(x,y)=1+o(1)\qquad \big((\log x)/y+(\log y)/\log x\to 0\big),$$

la même estimation étant valable pour $C^*(x,y)$.

La formule (1·13) traduit, sous une forme forte l'indépendance asymptotique des relations de divisibilité par les puissances de nombres premiers dans le domaine de friabilité indiqué. Cela induit naturellement la question de déterminer le domaine optimal pour ce phénomène.

Lorsque $u:=(\log x)/\log y$ est borné la valeur asymptotique de $C(x,y)$ a été déterminée dans [11] et tabulée dans [9]. On a alors

$$(1{\cdot}14)\qquad C(u):=\lim_{y\to\infty} C(y^u,y)\in\,]1,2],$$

et aussi $C^*(u):=\lim_{y\to\infty} C^*(y^u,y)=C(u)$.

Symétriquement, il est aisé de constater que, pour tout $y\geqslant 2$ fixé, nous avons

$$\liminf_{x\to\infty} C(x,y)\geqslant \mathrm{e}\Big(1-\frac{1}{\pi(y)}\Big)^{\pi(y)-1}>1,$$

avec la convention que le terme médian vaut e pour $y=2$. On obtient cette inégalité en choisissant $f(2^\nu)=1$ si $\nu=\nu(y):=\lfloor(\log x)/\{\pi(y)\log 2\}\rfloor$ et $f(p^j)=0$ si $p^j\neq 2^{\nu(y)}$.

En utilisant les estimations de [4], nous nous proposons ici de montrer le résultat suivant, qui illustre de manière particulièrement limpide, dans un domaine optimal, l'indépendance asymptotique des conditions de divisibilité par des puissances de nombres premiers distincts.

**Théorème 1.1.** *Sous les conditions $x\geqslant y\geqslant 2$ et $1/y+1/u=o(1)$, nous avons*

$$(1{\cdot}15)\qquad C(x,y)=1+o(1).$$

*De plus, la même relation est valable pour $C^*(x,y)$.*

Comme indiqué plus haut, seul le domaine complémentaire de celui de (1·13) est à considérer. La méthode employée ici permet en fait de conclure dès que $y\to\infty$, $y=o\big((\log x)^2/\log_2 x\big)$, et fournit pour le terme d'erreur de (1·15) l'estimation

$$\ll \frac{y\log y}{(\log x)^2}+\frac{\log y}{y}\cdot$$

Nous n'avons pas cherché dans cette note à expliciter le majorant du terme d'erreur de (1·15) produit par la concaténation des deux approches impliquées.

## 2. Preuve

Il résulte de ce qui précède que nous pouvons limiter l'étude au domaine $y \leqslant (\log x)^{3/2}$, dans lequel nous allons établir l'estimation

$$C(x,y) = 1 + O\Big(\frac{h}{\overline{u}^2}\Big) \tag{2·1}$$

où nous avons posé

$$h := \pi(y), \quad \overline{u} := \min(h,u).$$

Nous pouvons aussi, classiquement, nous restreindre au cas où $f$ est à valeurs réelles. Nous nous plaçons dans le cadre du modèle non biaisé, les détails essentiellement identiques dans le cas du modèle standard.

Avec les notations précédemment introduites, nous avons

$$\mathbb{E}(\xi_p^*) = \sum_{1\leqslant \nu\leqslant \nu_p} \frac{g_p(\alpha)f(p^\nu)}{p^{\nu\alpha}(1-w_p)}, \qquad \mathbb{V}(\xi_p^*) = g_p(\alpha)\mathbb{E}(\xi_p^*)^2 + \sum_{1\leqslant \nu\leqslant \nu_p} \frac{g_p(\alpha)F_p(\nu)^2}{p^{\nu\alpha}(1-w_p)},$$

où nous avons posé $F_p(\nu) := f(p^\nu) - \mathbb{E}(\xi_p^*)$ $(0 \leqslant \nu \leqslant \nu_p)$.

À ce stade, il est à noter, d'une part, que la fonction $F(n) := \sum_{p\leqslant y} F_p(v_p(n))$, où $v_p$ désigne l'application valuation $p$-adique, n'est en général pas additive sur $S(x,y)$ et, d'autre part, que, dans la décomposition canonique $f = g + h$ utilisée dans [5] avec $h$ fortement additive, la fonction $g$ vérifie

$$g(p^\nu) = f(p^\nu) - p^\alpha\mathbb{E}(\xi_p^*) = F_p(\nu) - g_p(\alpha)p^\alpha\mathbb{E}(\xi_p^*) \qquad (1 \leqslant \nu \leqslant \nu_p).$$

Alors que le modèle probabiliste de $g$ est d'espérance nulle, celui de $F$ est d'espérance proche de 0 — cf. (2·6) *infra.*

Cela étant, nous avons

$$\begin{aligned} V_f^* = V_f^*(x,y) &= \frac{1}{\Psi(x,y)} \sum_{n\in S(x,y)} \Big| \sum_{p\leqslant y} \Big\{ f\big(p^{v_p(n)}\big) - \mathbb{E}\big(\xi_p^*\big)\Big\}\Big|^2 \\ &= \sum_{p\leqslant y}\sum_{0\leqslant\nu\leqslant\nu_p} F_p(\nu)^2\frac{\Psi_p(x/p^\nu,y)}{\Psi(x,y)} + \sum_{\substack{p,q\leqslant y\\ p\neq q}} \sum_{\substack{0\leqslant\nu\leqslant\nu_p\\ 0\leqslant\mu\leqslant\nu_q}} F_p(\nu)F_q(\mu)\frac{\Psi_{pq}(x/p^\nu q^\mu,y)}{\Psi(x,y)} \\ &= A + B \quad \text{(disons)}, \end{aligned}$$

où, ici et dans la suite, nous notons traditionnellement $\Psi_m(x,y)$ le nombre des entiers $y$-friables qui sont premiers à un entier naturel $m$.

Le théorème 2.4 de [2] permet d'écrire

$$\frac{\Psi_p(x/p^\nu,y)}{\Psi(x,y)} = \big(1 - \tfrac12 t_p(\nu)^2\big)^{b\overline{u}}\frac{g_p(\alpha)}{p^{\nu\alpha}}\Big\{1 + O\Big(\frac{1}{\overline{u}} + t_p(\nu)\Big)\Big\} \qquad (x \geqslant y \geqslant 2),$$

avec la notation

$$t_p(\nu) := (\log p^\nu)/\log x$$

et où $b = b(p,\nu,x,y) \asymp 1$.

Comme $\{1+O(v)\}(1-\frac{1}{2}v^2)^{b\overline{u}} \leqslant 1+O\big(1/\overline{u}\big)$ pour $0 \leqslant v \leqslant 1$, il suit

$$A \leqslant \Big\{1+O\Big(\frac{1}{\overline{u}}\Big)\Big\}\sum_{p\leqslant y}\sum_{0\leqslant \nu\leqslant \nu_p}\frac{g_p(\alpha)F_p(\nu)^2}{p^{\nu\alpha}} \leqslant \Big\{1+O\Big(\frac{1}{\overline{u}}\Big)\Big\}\mathbb{V}(Z_f^*),$$

où nous avons fait usage de la majoration $w_p \ll \mathrm{e}^{-\overline{u}/2}$, établie dans [5], formules (1.12) et (1.13).

Pour évaluer $B$, nous faisons appel au lemme 2.1 de [4] que nous reformulons légèrement. Nous posons

$$v_k(\alpha) := \log k - \frac{g_k'(\alpha)}{g_k(\alpha)} = \log k - \sum_{p|k}\frac{\log p}{p^\alpha - 1} \qquad (k \geqslant 1),$$

$$(2{\cdot}2)\qquad \sigma_2 := \left|\left[\frac{\mathrm{d}}{\mathrm{d}s}\sum_{p\leqslant y}\frac{\log p}{p^s-1}\right]_{s=\alpha}\right| \asymp \frac{(u\log y)^2}{\overline{u}} \qquad (x \geqslant y \geqslant 2),$$

l'évaluation ayant été établie au lemme 4 de [10].

**Lemme 2.1.** *Sous les conditions* $x \geqslant y \geqslant 2$, $u > \sqrt{\log y}$, $p^\nu q^\mu \leqslant x$, $p \neq q \leqslant y$, *nous avons uniformément*

$$(2{\cdot}3)\qquad \begin{aligned}\frac{\Psi_{pq}\big(x/p^\nu q^\mu, y\big)}{\Psi(x,y)} &= \frac{g_{pq}(\alpha)}{p^{\nu\alpha}q^{\mu\alpha}}\Big\{1 - \frac{v_{p^\nu}(\alpha)v_{q^\mu}(\alpha)}{\sigma_2} + O(R)\Big\}\\ &\quad + \frac{g_p(\alpha)\Psi_q\big(x/q^\mu, y\big)}{p^{\nu\alpha}\Psi(x,y)} + \frac{g_q(\alpha)\Psi_p\big(x/p^\nu, y\big)}{q^{\mu\alpha}\Psi(x,y)},\end{aligned}$$

*où l'on a posé*

$$R := \frac{1}{\overline{u}^2} + \overline{u}^2(t_p(\nu)+t_q(\mu))^4 + \overline{u}^{5/2}(t_p(\nu)+t_q(\mu))^5.$$

Nous pouvons estimer la contribution à $B$ du domaine $p^\nu q^\mu > x^{1/10\sqrt{\overline{u}}}$ via la majoration

$$\frac{\Psi_{pq}(x/p^\nu q^\mu, y)}{\Psi(x,y)} \ll \frac{g_{pq}(\alpha)}{p^{\nu\alpha}q^{\mu\alpha}} \qquad (p^\nu q^\mu \leqslant x),$$

établie au théorème 2.4 de [3]. En observant que

$$\sum_{p^\nu > x^{1/20\sqrt{\overline{u}}}}\frac{g_p(\alpha)}{p^{\nu\alpha}} \ll \mathrm{e}^{-\sqrt{\overline{u}}/50}$$

nous obtenons que cette contribution est $\ll \mathbb{V}(Z_f^*)\mathrm{e}^{-\sqrt{\overline{u}}/60}$.

Afin d'évaluer la contribution à $B$ du terme d'erreur $R$ apparaissant dans (2·3), nous observons que, dans le domaine considéré,

$$(2{\cdot}4)\qquad y_p := \bigg(\sum_{0\leqslant \nu\leqslant \nu_p}\frac{g_p(\alpha)F_p(\nu)^2}{p^{\nu\alpha}}\bigg)^{1/2} \leqslant \sqrt{\mathbb{V}(\xi_p^*)},$$

$$(2{\cdot}5)\qquad z_p(j) := \bigg(\sum_{0\leqslant \nu\leqslant \nu_p}\frac{g_p(\alpha)t_p(\nu)^{2j}}{p^{\nu\alpha}}\bigg)^{1/2} \ll \frac{1}{\overline{u}^j} \qquad (0 \leqslant j \leqslant 5),$$

où la dernière majoration est issue du lemme 3.2 de [3]. Nous obtenons ainsi que cette contribution est

$$\begin{aligned}&\ll \frac{1}{\overline{u}^2}\sum_{p,q\leqslant y} y_py_qz_p(0)z_q(0)+\overline{u}^2\sum_{p,q\leqslant y} y_pz_p(4)y_qz_q(0)+\overline{u}^{5/2}\sum_{p,q\leqslant y} y_py_qz_p(5)z_q(0)\\ &\ll \frac{1}{\overline{u}^2}\Big(\sum_{p\leqslant y}\sqrt{\mathbb{V}(\xi_p^*)}\Big)^2\ll \frac{h\mathbb{V}(Z_f^*)}{\overline{u}^2}.\end{aligned}$$

Il reste à estimer la contribution à $B$ du terme principal de (2·3) lorsque $p^\nu q^\mu \leqslant x^{1/10\sqrt{\overline{u}}}$. À cette fin, nous réintroduisons les termes satisfaisant à la condition $p^\nu q^\mu > x^{1/10\sqrt{\overline{u}}}$, ce qui produit à nouveau une erreur acceptable.

Nous allons montrer que, à l'exclusion de celui impliquant $v_{p^\nu}(\alpha)v_{q^\mu}(\alpha)$, les termes principaux de (2·3) ne contribuent qu'au terme d'erreur de (2·1). L'identité

$$\sum_{0\leqslant\nu\leqslant\nu_p}\frac{F_p(\nu)g_p(\alpha)}{p^{\nu\alpha}}=\sum_{1\leqslant\nu\leqslant\nu_p}\frac{f(p^\nu)g_p(\alpha)}{p^{\nu\alpha}}-\mathbb{E}(\xi_p^*)\Big(1-\frac{w_p}{p^\alpha}\Big)=-g_p(\alpha)w_p\mathbb{E}(\xi_p^*) \tag{2·6}$$

nous permet en effet de majorer ces contributions par

$$\begin{aligned}&\ll \sum_{p\leqslant y}\sum_{\substack{q\leqslant y\\ q\neq p}}\sum_{0\leqslant\mu\leqslant\nu_q}\frac{g_q(\alpha)|F_q(\mu)|}{q^{\mu\alpha}}\Big|\sum_{0\leqslant\nu\leqslant\nu_p}\frac{g_p(\alpha)F_p(\nu)}{p^{\nu\alpha}}\Big|\\ &\ll \sum_{p\leqslant y}g_p(\alpha)w_p|\mathbb{E}(\xi_p^*)|\sum_{q\leqslant y}\Big(\sum_{0\leqslant\mu\leqslant\nu_q}\frac{g_q(\mu)F_q(\mu)^2}{q^{\mu\alpha}}\Big)^{1/2}\\ &\ll \mathrm{e}^{-\overline{u}/2}\Big(\sum_{p\leqslant y}g_p(\alpha)\Big)^{1/2}\Big(\sum_{p\leqslant y}g_p(\alpha)\mathbb{E}(\xi_p^*)^2\Big)^{1/2}\sum_{q\leqslant y}\mathbb{V}(\xi_q^*)^{1/2}\\ &\ll h\mathrm{e}^{-\overline{u}/2}\Big(\sum_{p\leqslant y}g_p(\alpha)\mathbb{E}(\xi_p^*)^2\Big)^{1/2}\mathbb{V}(Z_f^*)^{1/2}\ll h\mathrm{e}^{-\overline{u}/2}\mathbb{V}(Z_f^*)\ll \mathrm{e}^{-\overline{u}/3}\,\mathbb{V}(Z_f^*).\end{aligned}$$

Considérons enfin la quantité

$$C:=-\frac{1}{\sigma_2}\sum_{p\leqslant y}\sum_{\substack{q\leqslant p\\ q\neq p}}\sum_{0\leqslant\mu\leqslant\nu_q}\frac{g_q(\alpha)F_q(\mu)v_{q^\mu}(\alpha)}{q^{\mu\alpha}}\sum_{0\leqslant\nu\leqslant\nu_p}\frac{g_p(\alpha)F_p(\nu)v_{p^\nu}(\alpha)}{p^{\nu\alpha}}.$$

En réintroduisant les termes correspondant à $p=q$, et en majorant la contribution du carré ainsi obtenu par 0, nous obtenons

$$\begin{aligned}C&\leqslant \frac{1}{\sigma_2}\sum_{p\leqslant y}\Big(\sum_{0\leqslant\nu\leqslant\nu_p}\frac{g_p(\alpha)F_p(\nu)v_{p^\nu}(\alpha)}{p^{\nu\alpha}}\Big)^2\\ &\leqslant \frac{1}{\sigma_2}\sum_{p\leqslant y}\mathbb{V}(\xi_p^*)\sum_{0\leqslant\nu\leqslant\nu_p}\frac{g_p(\alpha)v_{p^\nu}(\alpha)^2}{p^{\nu\alpha}}\ll\frac{\mathbb{V}(Z_f^*)}{\overline{u}}.\end{aligned}$$

Nous avons ici utilisé (2·2), (2·5), et l'estimation $\alpha\gg\overline{u}/\log x$ $(x\geqslant y\geqslant 2)$, qui découle du lemme 3.1 de [3].

Nous avons ainsi établi que, pour une constante absolue convenable, nous avons

$$B\leqslant K\mathbb{V}(Z_f^*)h/\overline{u}^2$$

dans le domaine considéré.

Cela achève la démonstration de la majoration contenue dans (2·1).

Pour établir la minoration correspondante, nous considérons la fonction additive $f \in \mathbb{A}$ telle que $f(2) = 1$ et $f(p^\nu) = 0$ si $p = 2$, $\nu \geqslant 2$ ou $p \geqslant 3$, $\nu \geqslant 0$. Nous avons alors $B = 0$ et

$$\begin{aligned}
\mathbb{E}(\xi_2^*) &= \frac{g_2(\alpha)}{2^\alpha(1-w_2)},\\
A &= \mathbb{E}(\xi_2^*)^2\frac{\Psi_2(x,y)}{\Psi(x,y)} + \big\{1-\mathbb{E}(\xi_2^*)\big\}^2\frac{\Psi_2(x/2,y)}{\Psi(x,y)}\\
&= \mathbb{E}(\xi_2^*)^2 g_2(\alpha)\Big\{1+O\Big(\frac{1}{\overline{u}}\Big)\Big\} + \big\{1-\mathbb{E}(\xi_2^*)\big\}^2\frac{g_2(\alpha)}{2^\alpha}\Big\{1+O\Big(\frac{1}{\overline{u}}\Big)\Big\}\\
&= \mathbb{E}(\xi_2^*)^2 g_2(\alpha) + \big\{1-\mathbb{E}(\xi_2^*)\big\}^2\frac{g_2(\alpha)}{2^\alpha} + O\Big(\frac{g_2(\alpha)}{\overline{u}}\Big)
\end{aligned}$$

où l'avant-dernière estimation découle du théorème 2.4 de [3]. Comme

$$\mathbb{V}(Z_f^*) = \mathbb{E}(\xi_2^*)^2 g_2(\alpha) + \big\{1-\mathbb{E}(\xi_2^*)\big\}^2\frac{g_2(\alpha)}{2^\alpha(1-w_2)} \asymp g_2(\alpha)$$

et $w_2 \ll \mathrm{e}^{-\overline{u}/2}$, nous en déduisons que $V_f^*(x,y) = \mathbb{V}(Z_f^*)\big\{1+O\big(1/\overline{u}\big)\big\}$.

*Remarque.* Pour $j \geqslant 1$, le membre de droite de (2·5) peut être divisé par $p^{\alpha/2}$. Ce n'est pas le cas lorsque $j = 0$ puisque, par exemple, $g_p(\alpha) \asymp 1$ dès que $y > (\log x)^{1+\varepsilon}$ avec $\varepsilon > 0$ fixé. Comme $F_p(0)$ n'est en général pas nul, cela atténue, pour les grandes valeurs de $y$, la qualité du majorant de la contribution à $B$ du terme en $1/\overline{u}^2$ de $R$, pour finalement compromettre totalement la méthode lorsque $y/\log y \gg \overline{u}^2$.

**Remerciements.** Pendant les recherches ayant conduit au présent travail, le deuxième auteur a bénéficié d'une bourse du Conseil de Recherches en Sciences Naturelles et en Génie du Canada, et également d'un poste rouge du CNRS pour un séjour à l'institut Élie Cartan de Lorraine. Il tient à remercier cet établissement pour les exceptionnelles conditions de travail qui lui ont été offertes.

## Bibliographie

[1] K. Alladi, The Turán-Kubilius inequality for integers without large prime factors, *J. reine angew. Math.* **335** (1982), 180–196.

[2] R. de la Bretèche & G. Tenenbaum, Propriétés statistiques des entiers friables, *Ramanujan J.* **9** (2005), 139–202.

[3] R. de la Bretèche & G. Tenenbaum, Entiers friables : inégalité de Turán–Kubilius et applications, *Invent. Math.* **159** (2005), 531–588.

[4] R. de la Bretèche & G. Tenenbaum, On the friable Turán–Kubilius inequality, Proceedings of the Fifth International Conference *Analytic and Probabilistic Methods in Number Theory* held on the occasion of Professor Jonas Kubilius' 90th birthday (2012), E. Manstavičius et al. (eds), TEV Vilnius 2012, 259-265.

[5] R. de la Bretèche & G. Tenenbaum, Sur l'inégalité de Turán–Kubilius friable, *J. London Math. Soc.* (2) **93**, n° 1 (2016), 175–193.

[6] R. de la Bretèche & G. Tenenbaum, Une nouvelle approche dans la théorie des entiers friables, *Compositio Math.* **153** (2017), 453-473.

[7] P.D.T.A. Elliott, *Probabilistic number theory : mean value theorems*, Grundlehren der Math. Wiss. 239, Springer-Verlag, New York, Berlin, Heidelberg 1979.

[8] P.D.T.A. Elliott, *Probabilistic number theory : central limit theorems*, Grundlehren der Math. Wiss. 240, Springer-Verlag, New York, Berlin, Heidelberg 1980.

[9] G. Hanrot, B. Martin & G. Tenenbaum, Constantes de Turán-Kubilius friables : étude numérique, *Exp. Math.* Experiment. Math. **19, n° 3** (2010), 345–361.

[10] A. Hildebrand & G. Tenenbaum, On integers free of large prime factors, *Trans. Amer. Math. Soc.* **296** (1986), 265–290.

[11] B. Martin & G. Tenenbaum, Sur l'inégalité de Turán-Kubilius friable, *J. reine angew. Math.* **647** (2010), 175–234.

[12] G. Tenenbaum, *Introduction à la théorie analytique et probabiliste des nombres*, quatrième édition, coll. Échelles, Belin, 2015, 592 pp

[13] T.Z. Xuan, The Turán-Kubilius inequality for integers free of large prime factors, *J. Number Theory* **43** (1993), 82–87.

[14] T.Z. Xuan, The Turán-Kubilius inequality for integers free of large prime factors (II), *Acta Arith.* **65** (1993), 329–352.

Régis de la Bretèche
Institut de Mathématiques
de Jussieu-PRG
UMR 7586
Université Paris Diderot-Paris 7
Sorbonne Paris Cité,
Case 7012, F-75013 Paris
France

Youness Lamzouri
Department of Mathematics
and Statistics
York University
4700 Keele St
Toronto, ON, M3J1P3,
Canada

Gérald Tenenbaum
Institut Élie Cartan
Université de Lorraine
BP 70239
54506 Vandœuvre-lès-Nancy Cedex
France